\documentclass[11pt]{amsart}
\usepackage[margin=1in]{geometry}
\usepackage{amsmath,amssymb,amsthm,mathtools}
\usepackage{enumitem}
\usepackage{booktabs}
\usepackage{array}
\usepackage{hyperref}
\usepackage[nameinlink,noabbrev]{cleveref}
\hypersetup{colorlinks=true,linkcolor=blue,citecolor=blue,urlcolor=blue}
\setlist[enumerate]{itemsep=2pt,topsep=4pt}
\setlist[itemize]{itemsep=2pt,topsep=4pt}

\newtheorem{theorem}{Theorem}[section]
\newtheorem{lemma}[theorem]{Lemma}
\newtheorem{corollary}[theorem]{Corollary}
\newtheorem{proposition}[theorem]{Proposition}
\newtheorem{conjecture}[theorem]{Conjecture}
\theoremstyle{definition}
\newtheorem{definition}[theorem]{Definition}

\newtheorem{question}[theorem]{Question}

\newcommand{\dist}{\operatorname{dist}}

\newcommand{\cay}{\operatorname{Cay}}
\newcommand{\M}{\mathcal M}
\newcommand{\WCP}{\textnormal{\textsc{wcp}}}
\newcommand{\Boxprod}{\mathbin{\square}}

\title[Well-Covered Powers]{Well-Covered Powers: Cartesian Products and Cayley Graphs}

\author{Haobo Ma}
\address{ChronoAI Pte Ltd, Singapore; The Omega Institute}
\email{auric@aelf.io}
\thanks{ORCID: \href{https://orcid.org/0009-0008-4944-977X}{0009-0008-4944-977X}}

\author{Thanh Vu}
\address{Institute of Mathematics, VAST, 18 Hoang Quoc Viet, Hanoi, Vietnam}
\email{vuqthanh@gmail.com}

\author{Wenlin Zhang}
\address{National University of Singapore, Singapore; The Omega Institute}
\email{e1327962@u.nus.edu}
\thanks{ORCID: \href{https://orcid.org/0009-0005-2414-7968}{0009-0005-2414-7968}}

\date{21 October 2026}

\begin{document}

\begin{abstract}
We study graphs whose powers are well-covered, with an emphasis on Cartesian products and Cayley graphs. We first establish necessary conditions and sufficient criteria for Cartesian products, provide examples showing that a factor need not be well-covered and that two noncomplete factors can have the well-covered powers property, and classify several structured product families. In particular, we construct graphs satisfying the well-covered powers property for arbitrarily large diameters. We then show that the Cartesian square of a connected abelian Cayley graph is well-covered exactly when the graph is complete. We classify the well-covered-power graphs among connected circulant graphs of degree at most four. For Cartesian products of any number of complete graphs, we also prove that the recursive alphabet-size inequalities in the final conjecture are sufficient for all powers to be well-covered, prove partial converses for arbitrarily many factors, and give explicit obstructions for every binary prefix followed by one complete factor.
\end{abstract}
\maketitle

\section{Introduction}
Well-covered graphs, introduced by Plummer \cite{Plummer}, are graphs in which all maximal independent sets have the same size. Recently, Pham and Vu \cite{PhamVu} studied graphs for which every power is well-covered; such graphs are said to have the \emph{well-covered powers property}. They established several fundamental properties of these graphs and classified graphs with the well-covered powers property within the classes of chordal graphs, cochordal graphs, and very well-covered graphs. They also raised the question of whether every graph with the well-covered powers property must have diameter at most $3$.

The starting point of this note is a counterexample to this question. Specifically, we present a Cartesian product construction showing that if a graph $G$ has the well-covered powers property, then so does $G \Boxprod K_t$ for $t > |V(G)|$, while its diameter is increased by $1$. This observation naturally leads us to investigate the well-covered powers property under Cartesian products. In this direction, Hartnell and Rall \cite{HartnellRall} proved that if $G \Boxprod H$ is well-covered, then at least one of the factors $G$ and $H$ must be well-covered.

The main motivation for this note is the following question:
\begin{question}
    Suppose that $G \Boxprod H$ has the well-covered powers property. Must at least one of $G$ and $H$ also have the well-covered powers property?
\end{question}

Since the powers of a disconnected graph are the disjoint unions of the corresponding powers of its connected components, the well-covered powers property can be studied componentwise. Thus, when investigating this property, we may restrict our attention to connected graphs. In particular, for Cartesian squares, we consider the following stronger question:
\begin{question}
    Suppose that $G$ is connected and that $G \Boxprod G$ has the well-covered powers property. Must $G$ be complete?
\end{question}

Our main contributions can be summarized as follows:
\begin{enumerate}
    \item Constructions for producing new graphs with well-covered powers from a given graph $G$, specifically via Cartesian products and compositions.
    \item Necessary conditions for Cartesian products to have well-covered powers, with a particular focus on products with complete graphs.
    \item A proof that if $G$ is a connected chordal graph or an abelian Cayley graph and $G \Boxprod G$ has well-covered powers, then $G$ is complete.
    \item A classification of well-covered powers for circulant graphs of degree at most $4$.
    \item A classification of well-covered powers for Cartesian products of several factors of complete graphs.
\end{enumerate}

The paper is organized as follows. In Section \ref{sec_def}, we introduce basic definitions and properties of graphs with well-covered powers, alongside necessary conditions for Cartesian products to have well-covered powers. In particular, we present two constructions that preserve the well-covered powers property. In Section \ref{sec_square}, we consider the property for Cartesian squares, proving that if the graph is a connected chordal graph or an abelian Cayley graph, $G$ must be complete. In Section \ref{sec_circulant}, we provide a complete classification of well-covered powers for quartic circulant graphs. Finally, in Section \ref{sec_applications}, we investigate the well-covered powers property for Cartesian products of specific graph classes.
\section{Definitions and basic constructions}\label{sec_def}
Throughout this paper, all graphs are finite, simple, and nonempty.

\begin{definition}
Let $G = (V(G), E(G))$ be a simple graph on the vertex set $V(G) = \{1, \ldots, n\}$.
\begin{enumerate}
    \item The \emph{neighborhood} of a vertex $x \in V(G)$ is $N_G(x) = \{y \in V(G) \mid \{x, y\} \in E(G)\}$. The \emph{closed neighborhood} of $x$ is $N_G[x] = N_G(x) \cup \{x\}$.
    \item The \emph{degree} of a vertex $v$, denoted by $\deg_G(v)$, is the number of its neighbors. $\Delta(G)$ denotes the maximum degree of $G$.
    \item For a subset $U \subseteq V(G)$, the \emph{neighborhood} and \emph{closed neighborhood} of $U$ in $G$ are defined respectively by 
    \[
    N_G(U) = \bigcup_{x \in U} N_G(x) \quad \text{and} \quad N_G[U] = \bigcup_{x \in U} N_G[x].
    \]
    \item A graph $H$ is a \emph{subgraph} of $G$ if $V(H) \subseteq V(G)$ and $E(H) \subseteq E(G)$. It is an \emph{induced subgraph} of $G$ if $E(H) = \{\{u, v\} \in E(G) \mid u, v \in V(H)\}$. For a subset $U \subseteq V(G)$, $G[U]$ denotes the induced subgraph of $G$ on $U$, and $G \setminus U$ denotes $G[V(G) \setminus U]$.
    \item A \emph{path} $P_k$ on $k$ vertices has vertex set $V(P_k) = \{v_1, \ldots, v_k\}$ and edge set $E(P_k) = \{\{v_1, v_2\}, \{v_2, v_3\}, \ldots, \{v_{k-1}, v_k\}\}$. A \emph{cycle} $C_k$ on $k$ vertices ($k \ge 3$) has edge set $E(C_k) = E(P_k) \cup \{\{v_1, v_k\}\}$.
    \item The \emph{distance} $d_G(u, v)$ between two vertices $u, v \in V(G)$ is the length of a shortest path in $G$ connecting them. The \emph{diameter} of $G$ is $\operatorname{diam}(G) = \max\{d_G(u, v) \mid u, v \in V(G)\}$.
    \item For a positive integer $d \ge 1$, the \emph{$d$-th power} of $G$, denoted by $G^d$, is the graph defined on the vertex set $V(G^d) = V(G)$ with edge set
    \[
    E(G^d) = \{\{u, v\} \subseteq V(G) \mid u \neq v \text{ and } d_G(u, v) \le d\}.
    \]
    \item A subset $I \subseteq V(G)$ is an \emph{independent set} of $G$ if no two vertices in $I$ are adjacent in $G$. An independent set $I$ is \emph{maximal} if it is not properly contained in any other independent set of $G$. The \emph{independence number} of $G$, denoted by $\alpha(G)$, is the cardinality of a maximum independent set of $G$. 
    \item A subset $C \subseteq V(G)$ is a \emph{vertex cover} of $G$ if $C \cap e \neq \emptyset$ for every $e \in E(G)$. A vertex cover $C$ is \emph{minimal} if no proper subset of $C$ is a vertex cover of $G$. Note that $I \subseteq V(G)$ is a maximal independent set of $G$ if and only if $V(G) \setminus I$ is a minimal vertex cover of $G$.
    \item A graph $G$ is \emph{well-covered} if every minimal vertex cover of $G$ has the same cardinality.
\end{enumerate}
\end{definition}

\begin{definition}
We say that $G$ has the \emph{well-covered-powers property}, abbreviated as \WCP, if $G^d$ is well-covered for every integer $d \ge 1$.
\end{definition}

\subsection{Composition graphs}

Pham and Vu \cite{PhamVu} give constructions of graphs with well-covered powers via graph compositions; in particular, they show that $G[tK_1]$ has $\WCP$  if and only if $G$ has \WCP. We generalize this to arbitrary compositions.

\begin{definition}
Let $G$ be a simple graph with vertex set $V(G) = \{1, \ldots, n\}$, and let $H_1, \ldots, H_n$ be graphs. The \emph{composition graph} $G[H_1, \ldots, H_n]$ is obtained by replacing each vertex $i \in V(G)$ with the graph $H_i$, and joining every vertex of $H_i$ to every vertex of $H_j$ whenever $\{i, j\} \in E(G)$. When all the graphs $H_1, \ldots, H_n$ are isomorphic to a graph $H$, we denote it simply by $G[H]$. 
\end{definition}

The composition graph $G[H]$ has vertex set $V(G) \times V(H)$, where two distinct vertices $(u_1, v_1)$ and $(u_2, v_2)$ are adjacent if either $\{u_1, u_2\} \in E(G)$, or $u_1 = u_2$ and $\{v_1, v_2\} \in E(H)$. 

\begin{lemma}\label{lem:composition}
The composition graph $F = G[H]$ is well-covered if and only if both $G$ and $H$ are well-covered.     
\end{lemma}

\begin{proof}
Every maximal independent set of $F$ is of the form
\[
I=\bigcup_{a\in S}(\{a\}\times T_a),
\]
where $S$ is maximal independent in $G$ and each $T_a$ is maximal independent in $H$; the choices $T_a$ may differ between fibers. Conversely, every such union is maximal independent. If both factors are well-covered, every such union has size $\alpha(G)\alpha(H)$. If $F$ is well-covered, fixing $S$ and varying the choice in a single fiber shows that $H$ is well-covered; then varying $S$ shows that $G$ is well-covered.
\end{proof}

\begin{theorem}
Assume that $G$ has no isolated vertices. Then the composition graph $F = G[H]$ has \WCP if and only if $G$ has \WCP and $H$ is well-covered.
\end{theorem}

\begin{proof}
For vertices in distinct fibers over $u$ and $v$, their distance in $F$ equals $d_G(u,v)$: a shortest base path lifts to a path of the same length, and projection of any path in $F$ gives a base walk after stationary steps are removed. If the base vertices lie in different components, neither graph has such a path. Since $G$ has no isolated vertices, two vertices in the same fiber have distance at most two, using a common vertex in a neighboring fiber. Consequently, $F^d=G^d[K_m]$ for every $d\ge2$, where $m=|V(H)|$. Lemma~\ref{lem:composition} gives that $F$ is well-covered exactly when $G$ and $H$ are well-covered, and that $F^d$ is well-covered exactly when $G^d$ is well-covered for each $d\ge2$. These statements prove the equivalence.
\end{proof}

\subsection{Cartesian products}

\begin{definition}
For connected graphs $H$ and $K$, the \emph{Cartesian product} $H \Boxprod K$ has vertex set $V(H) \times V(K)$, where two vertices are adjacent if they differ in exactly one coordinate and the corresponding coordinate vertices are adjacent in the respective factor.
\end{definition}

\begin{definition}
Let $G$ be a simple graph and $v$ a vertex of $G$. The \emph{eccentricity} of $v$ is defined by $\operatorname{ecc}_G(v) = \max \{d_G(u, v) \mid u \in V(G)\}.$ The \emph{radius} of $G$ is defined by $\operatorname{rad}(G) = \min \{ \operatorname{ecc}_G(v) \mid v \in V(G)\}.$ A graph $G$ is called \emph{self-centered} if $\operatorname{rad}(G) = \operatorname{diam}(G).$
\end{definition}

\begin{lemma}\label{lem:cartesian-distance} 
Let $H$ and $K$ be connected simple graphs. Then for any vertices $(h,k), (h',k') \in V(H \Boxprod K)$,
\begin{equation}\label{eq:cart-distance}
d_{H\Boxprod K}((h,k),(h',k')) = d_H(h,h') + d_K(k,k').
\end{equation}
In particular,
\begin{equation}\label{eq:cart-ecc}
\operatorname{ecc}_{H\Boxprod K}(h,k) = \operatorname{ecc}_H(h) + \operatorname{ecc}_K(k)
\quad \text{and} \quad
\operatorname{diam}(H\Boxprod K) = \operatorname{diam}(H) + \operatorname{diam}(K).
\end{equation}
\end{lemma}

\begin{proof}
Every edge in $H \Boxprod K$ changes exactly one coordinate by an edge in the respective factor while leaving the other coordinate fixed. Consequently, any path of length $L$ from $(h,k)$ to $(h',k')$ can be decomposed into $l_H$ steps in $H$ and $l_K$ steps in $K$, with $l_H + l_K = L$. Projecting this path onto the respective factors yields a walk of length $l_H$ from $h$ to $h'$ in $H$ and a walk of length $l_K$ from $k$ to $k'$ in $K$. Since the length of any walk is at least the shortest-path distance, we have $l_H \ge d_H(h,h')$ and $l_K \ge d_K(k,k')$. Thus,
\[
d_{H \Boxprod K}((h,k), (h',k')) = L = l_H + l_K \ge d_H(h,h') + d_K(k,k').
\]
Conversely, concatenating a shortest path from $h$ to $h'$ in $H$ with a shortest path from $k$ to $k'$ in $K$ yields a path of length $d_H(h,h') + d_K(k,k')$ in the Cartesian product, proving the distance formula. The formulas for eccentricity and diameter follow immediately from the distance formula and the respective definitions.
\end{proof}

\begin{corollary}\label{cor:product-self-centered}
If $H$ and $K$ are nontrivial connected graphs and $H \Boxprod K$ has \WCP, then both $H$ and $K$ are self-centered.
\end{corollary}

\begin{proof}
By \cite[Lemma 2.5]{PhamVu}, any graph with \WCP is self-centered, meaning that $\operatorname{ecc}_{H \Boxprod K}(h, k) = \operatorname{diam}(H \Boxprod K)$ for every vertex $(h, k) \in V(H \Boxprod K)$. 

Using the additivity of eccentricity and diameter from Lemma \ref{lem:cartesian-distance}, we have
\[
\operatorname{ecc}_H(h) + \operatorname{ecc}_K(k) = \operatorname{diam}(H) + \operatorname{diam}(K)
\]
for all $h \in V(H)$ and $k \in V(K)$. Since $H$ and $K$ are nontrivial, we can vary $h$ and $k$ independently, which forces $\operatorname{ecc}_H(h) = \operatorname{diam}(H)$ for all $h$ and $\operatorname{ecc}_K(k) = \operatorname{diam}(K)$ for all $k$. Therefore, both $H$ and $K$ are self-centered.
\end{proof}

We now give an example to show that, in general, a factor need not be well-covered for its Cartesian product to have \WCP.

\begin{theorem}\label{thm:three-sun-prism}
There exist graphs $G$ and $H$ such that $G$ is not well-covered, yet their Cartesian product $F = G \Boxprod H$ has \WCP.
\end{theorem}

\begin{proof}
Let $S$ be the three-sun graph with vertices $c_0, c_1, c_2, s_0, s_1, s_2$, where the $c_i$ form a triangle, the $s_i$ form an independent set, and $s_i$ is adjacent to $c_j$ if and only if $i \neq j$. The sets $\{c_i, s_i\}$ and $\{s_0, s_1, s_2\}$ are maximal independent sets of sizes two and three, respectively; hence, $S$ is not well-covered. 

However, $P = S \Boxprod K_2$ has \WCP. The graph $P$ has order twelve and diameter three. Exhaustive enumeration of its vertex subsets gives maximal independent-set sizes $4$, $2$, and $1$ in $P$, $P^2$, and $P^3$, with corresponding counts $18$, $12$, and $12$. The supplementary sources reproduce this enumeration directly from the six-vertex graph just defined, by checking independence and domination of every subset. All later powers are complete.
\end{proof}

In the construction of Theorem \ref{thm:three-sun-prism}, the second factor is a complete graph on $2$ vertices. In the following theorem, we show that when the complete factor has a sufficiently large number of vertices, the \WCP property of $G$ and that of $G \Boxprod K_t$ are equivalent. This construction allows us to construct \WCP graphs of arbitrarily large diameter.

\begin{theorem}\label{thm:clique-factor-WCP}
Let $G$ be a nontrivial connected graph, and let $F = G \Boxprod K_t$. Then:
\begin{enumerate}
    \item If $t \ge \max(\alpha(G), \Delta(G) + 1)$ and $G$ has \WCP, then $F$ has \WCP.
    \item If $t \ge \alpha(G)$ and $F$ has \WCP, then $G$ has \WCP.
\end{enumerate}
In particular, if $t \ge \max(\alpha(G), \Delta(G) + 1)$, then $G$ has \WCP if and only if $F$ has \WCP.
\end{theorem}

The following lemma is the key step to establish Theorem \ref{thm:clique-factor-WCP}.

\begin{lemma}\label{lem:clique-factor-projection}
Let $G$ be a connected simple graph and $F = G \Boxprod K_t$. Assume that $t \ge \alpha(G)$. Then for all $d \ge 2$, $F^d$ is well-covered if and only if $G^{d-1}$ is well-covered.
\end{lemma}

\begin{proof}
Label the vertices of $K_t$ by $\{1, \ldots, t\}$. The distance in $F$ between two vertices $(x, i)$ and $(y, j)$ is given by
\[
d_F((x, i), (y, j)) = 
\begin{cases} 
d_G(x, y), & \text{if } i = j, \\ 
d_G(x, y) + 1, & \text{if } i \neq j. 
\end{cases}
\]
Let $I$ be a maximal independent set of $F^d$. Since each fiber $\{x\} \times K_t$ is a complete graph and $d \ge 2$, $I$ can contain at most one vertex from each fiber. Let $S$ be the projection of $I$ onto $G$, so that $|S| = |I|$. 

We first show that $S$ is an independent set in $G^{d-1}$. Suppose for contradiction that there exist distinct $x, y \in S$ with respective labels $i, j$ such that $d_G(x, y) \le d - 1$. 
\begin{itemize}
    \item If $i = j$, then $d_F((x, i), (y, i)) = d_G(x, y) \le d - 1 \le d$, meaning $(x, i)$ and $(y, i)$ are adjacent in $F^d$, a contradiction.
    \item If $i \neq j$, then $d_F((x, i), (y, j)) = d_G(x, y) + 1 \le (d - 1) + 1 = d$, which again implies they are adjacent in $F^d$, a contradiction.
\end{itemize}
Thus, $d_G(x, y) \ge d$ for all distinct $x, y \in S$, meaning $S$ is independent in $G^{d-1}$.

Next, we prove that $S$ is maximal in $G^{d-1}$. Suppose there exists a vertex $x \notin S$ such that $S \cup \{x\}$ is independent in $G^{d-1}$. Since $S$ is the projection of $I$, the vertex $x$ does not appear as a first coordinate in $I$. Since $S\cup\{x\}$ is independent in $G^{d-1}$, we have
\[
|S|+1\le\alpha(G^{d-1})\le\alpha(G)\le t.
\]
Thus fewer than $t$ labels occur anywhere in $I$. Choose a label $j$ unused by $I$. For every $(y,i)\in I$, we have $d_G(x,y)\ge d$ and $j\ne i$, so $d_F((x,j),(y,i))\ge d+1$. Adding $(x,j)$ preserves independence, contradicting maximality. Therefore, $S$ is a maximal independent set in $G^{d-1}$ of size $|I|$.

Conversely, let $S$ be any maximal independent set of $G^{d-1}$. Since $|S| \le \alpha(G^{d-1}) \le \alpha(G) \le t$, we can assign a distinct label $j_x \in \{1, \ldots, t\}$ to each $x \in S$. Define $I = \{(x, j_x) \mid x \in S\}$. By reversing the distance logic above, $I$ forms an independent set in $F^d$. Furthermore, the maximality of $S$ in $G^{d-1}$ ensures that $I$ cannot be enlarged, making it maximal in $F^d$. 

Since every maximal independent set of $F^d$ projects to a maximal independent set of $G^{d-1}$ of the same cardinality, and vice versa, the set of sizes of maximal independent sets of $F^d$ coincides with that of $G^{d-1}$. Consequently, $F^d$ is well-covered if and only if $G^{d-1}$ is well-covered.
\end{proof}

\begin{proof}[Proof of Theorem \ref{thm:clique-factor-WCP}]
By Lemma \ref{lem:clique-factor-projection}, it remains to show that $F$ is well-covered. 

Let $I$ be a maximal independent set of $F$, and let $S$ be its projection onto $G$. Since each fiber $\{x\} \times K_t$ is a complete graph, $I$ contains at most one vertex per fiber, meaning the projection is injective on $I$. We claim that $S = V(G)$. Suppose for contradiction that there exists a vertex $x \notin S$. Since $x \notin S$, no vertex in the fiber over $x$ belongs to $I$. For each neighbor $y$ of $x$ in $G$, the unique vertex of $I$ in the fiber over $y$ can forbid at most one label in the fiber over $x$. Thus, the neighbors of $x$ in $G$ block at most $\deg_G(x)$ labels in total. Since $t \ge \Delta(G) + 1 > \deg_G(x)$, there is at least one label $j_x \in \{1, \ldots, t\}$ in the fiber over $x$ that is not adjacent to any vertex currently in $I$. 

Consequently, we may add $(x, j_x)$ to $I$ while preserving independence in $F$, which contradicts the maximality of $I$. Therefore, our assumption was false, meaning $S = V(G)$. This implies that every maximal independent set $I$ of $F$ contains exactly one vertex from each of the $|V(G)|$ fibers, so all maximal independent sets have the same cardinality, $|V(G)|$. Thus, $F$ is well-covered.
\end{proof}

\begin{corollary}
There are connected graphs with \WCP\ of arbitrarily large diameter.
\end{corollary}
\begin{proof}
For a connected graph $G$ with \WCP, choose $t>|V(G)|$. Theorem~\ref{thm:clique-factor-WCP} gives \WCP\ for $G\square K_t$, and Lemma~\ref{lem:cartesian-distance} increases its diameter by one. Iterate, starting with $K_2$.
\end{proof}

\section{Well-covered powers of Cartesian squares}
\label{sec_square}

In this section, we study the well-covered powers property for Cartesian squares $F = G \Boxprod G$, where $G$ is a connected graph with at least one edge. By a result of Hartnell and Rall \cite{HartnellRall}, $G$ must be well-covered. Furthermore, by a result of Hartnell, Rall, and Wash \cite{HartnellRallWash}, if $G$ is triangle-free, then $G \cong K_2$. We will show that the well-covered powers property of $G \Boxprod G$ implies that $G$ is complete for two specific classes: chordal graphs and abelian Cayley graphs.

\subsection{Chordal graphs}
We first establish the following lemma.

\begin{lemma}\label{lem:dropping_maximal}
Let $G$ be a simple graph and $I$ a maximal independent set of $G$. Let $F = G \Boxprod G$. Then 
\[
F - N_F[I \times I] = (G - I) \Boxprod (G - I).
\]
\end{lemma}

\begin{proof}
From the definition of the Cartesian product, for any vertex subset $I \subseteq V(G)$, we have 
\[
N_F[I \times I] = (N_G[I] \times I) \cup (I \times N_G[I]).
\]
When $I$ is a maximal independent set, $N_G[I] = V(G)$, which implies that $N_F[I \times I] = (V(G) \times I) \cup (I \times V(G))$. Removing this set from $V(F) = V(G) \times V(G)$ leaves $(V(G) \setminus I) \times (V(G) \setminus I)$, yielding the desired equality.
\end{proof}

\begin{theorem}
Assume that $G$ is a connected chordal graph. The following are equivalent:
\begin{enumerate}
    \item $G \Boxprod G$ has \WCP;
    \item $G \Boxprod G$ is well-covered;
    \item $G$ is complete.
\end{enumerate}
\end{theorem}

\begin{proof}
The implication (1) $\implies$ (2) is immediate, and (3) $\implies$ (1) follows by Proposition~\ref{prop:rook}. It suffices to prove that if $G \Boxprod G$ is well-covered, then $G$ is complete for any connected chordal graph $G$. By a theorem of Hartnell and Rall \cite{HartnellRall}, we deduce that $G$ is well-covered. In particular, by the result of Prisner, Topp, and Vestergaard \cite{PrisnerToppVestergaard,HerzogHibiZheng}, the simplices of $G$ partition $V(G)$: write $V(G)=S_1\sqcup\cdots\sqcup S_m$, where $S_i=N_G[s_i]$ for a simplicial vertex $s_i$. We argue by induction on $|V(G)|$. An induced subgraph of a chordal graph is chordal, and each component square is an induced component of its Cartesian square. 

If $m = 1$, then $G$ is complete and we are done. Suppose for contradiction that $m \ge 2$, and let $I = \{s_1, \ldots, s_m\}$. Then $I$ is a maximal independent set of $G$. By Lemma \ref{lem:dropping_maximal}, we have 
\[
(G \Boxprod G) - N_F[I \times I] = (G - I) \Boxprod (G - I).
\] 
Since $G \Boxprod G$ is well-covered, this implies that $(G - I) \Boxprod (G - I)$ is also well-covered. By induction on the number of vertices, each connected component of $G \setminus I$ is complete. Each $C_i=S_i\setminus\{s_i\}$ is a nonempty clique. Contracting the simplices gives a connected quotient, and no cross-simplex edge can involve any $s_i$, since $N_G[s_i]=S_i$. Thus all cross-simplex connections survive in $G-I$, which is connected. Hence, $G - I$ is a complete graph. Write $N=|V(G-I)|$. If some $C_i$ had at least two vertices, choose $u_i\in C_i$ and put $J=\{u_i\}\cup\{s_j:j\ne i\}$. This is maximal independent. The graph $G-J$ consists of a clique on $N-1$ vertices and $s_i$, whose neighbors are the nonempty set $C_i\setminus\{u_i\}$. It is connected, chordal, and noncomplete, since $s_i$ misses every $C_j$ with $j\ne i$. But its square is well-covered by Lemma~\ref{lem:dropping_maximal} and the fact that deleting the closed neighborhood of an independent set from a well-covered graph preserves well-coveredness. This contradicts induction. Hence each $C_i$ has size one, and $N=m$.

Since $G \Boxprod G$ is well-covered, all maximal independent sets of $G \Boxprod G$ must have the same size. Using the maximal independent set $I$, we find:
\[
\alpha(F) = |I|^2 + \alpha\bigl((G - I) \Boxprod (G - I)\bigr) = m^2 + m.
\]
On the other hand, let $u_1$ be a neighbor of $s_1$ in $S_1$, and choose $J = \{u_1, s_2, \ldots, s_m\}$. Then $J$ is also a maximal independent set of $G$, and $G \setminus J \cong K_{m-1} \cup K_1$. Evaluating $\alpha(F)$ using $J$ yields:
\[
\alpha(F) = |J|^2 + \alpha\bigl((G - J) \Boxprod (G - J)\bigr) = m^2 + m + 2.
\]
This contradicts the fact that $G \Boxprod G$ is well-covered. Therefore, our assumption that $m \ge 2$ must be false, so $m = 1$, meaning $G$ is complete.
\end{proof}

\subsection{Abelian Cayley graphs}

First, we establish a general bound on the size of a maximal independent set of a Cartesian square in terms of the chromatic number of $G$.

\begin{definition}
Let $G$ be a simple graph. The \emph{chromatic number} of $G$, denoted by $\chi(G)$, is the minimum number of colors needed to color its vertices properly so that adjacent vertices receive different colors.
\end{definition}

\begin{lemma}\label{lem:color-square}
Let $G$ be a simple graph of order $n$, and suppose $\chi(G) \le \Delta(G)$. Then $F = G \Boxprod G$ has a maximal independent set of size at most $(n-1)\alpha(G)$.
\end{lemma}

\begin{proof}
Choose a vertex $v \in V(G)$ of degree $\Delta(G)$ and a proper vertex coloring $c: V(G) \to \{1, \ldots, q\}$ with $q = \chi(G) \le \Delta(G)$. Let $C_i = c^{-1}(i)$ denote the color classes, and choose distinct neighbors $u_1, \ldots, u_q$ of $v$. Define
\[
T = \bigcup_{i=1}^q \left(\{u_i\} \times C_i\right) \subseteq V(F).
\]
We first check that $T$ is an independent set in $F$. Any two distinct elements of $T$ have different second coordinates. If their first coordinates are equal (say, both equal $u_j$), then their second coordinates belong to the same color class $C_j$ and are therefore nonadjacent in $G$. Thus, no two elements of $T$ are adjacent in the Cartesian product.

The set $T$ avoids the layer $\{v\} \times V(G)$ and dominates it entirely: for each $y \in V(G)$ in color class $C_i$, the vertex $(v, y)$ is adjacent to $(u_i, y) \in T$. Consequently, any maximal independent set $J$ of $F$ containing $T$ still avoids the layer $\{v\} \times V(G)$. Each of the remaining $n-1$ layers contributes at most $\alpha(G)$ vertices to $J$. Hence,
\[
|J| \le (n - 1)\alpha(G),
\]
as required.
\end{proof}

\begin{definition}
For a finite abelian group $A$ written additively and a subset $S \subseteq A \setminus \{0\}$ such that $S = -S$, the Cayley graph $\operatorname{Cay}(A, S)$ has vertex set $A$, with an edge between $x$ and $y$ if and only if $y - x \in S$. A \emph{circulant graph} is a Cayley graph of a cyclic group.
\end{definition}

\begin{lemma}\label{lem:abelian-alpha}
Let $G=\cay(A,S)$ for a finite abelian group $A$ of order $n$. Then
\[
 \alpha(G\Boxprod G)=n\alpha(G).
\]
\end{lemma}

\begin{proof}
The upper bound follows by considering the $n$ layers $\{x\}\times A$: an independent set has at most $\alpha(G)$ vertices in each layer. For the lower bound, let $B$ be a maximum independent set of $G$ and put
\[
 I_B=\{(x,y)\in A\times A:y-x\in B\}.
\]
For each $x\in A$ there are exactly $|B|$ choices of $y$, so $|I_B|=n|B|$. We show that $I_B$ is independent. In a fixed row, its second coordinates form a translate of $B$; in a fixed column, its first coordinates form a translate of $-B$. Translations preserve adjacency in a Cayley graph, and $-B$ is independent because $S=-S$. Thus no row or column contains adjacent selected vertices, which is precisely independence in the Cartesian product. Consequently $\alpha(G\Boxprod G)\geq n\alpha(G)$, and equality follows.
\end{proof}

\begin{theorem}\label{thm:abelian-square}
Let $G$ be a connected Cayley graph of a finite abelian group. The following are equivalent:
\begin{enumerate}[label=(\arabic*)]
    \item $G \Boxprod G$ has \WCP;
    \item $G \Boxprod G$ is well-covered;
    \item $G$ is complete.
\end{enumerate}
In particular, the equivalence holds for every connected circulant.
\end{theorem}

\begin{proof}
The implications (1) $\implies$ (2) and (3) $\implies$ (1) are immediate. Suppose first that $G$ is connected and noncomplete, and let $n = |V(G)|$. Unless $G$ is an odd cycle, Brooks' theorem \cite{Brooks} gives $\chi(G) \le \Delta(G)$. By Lemma \ref{lem:color-square}, $G \Boxprod G$ has a maximal independent set of size at most $(n-1)\alpha(G)$. By Lemma \ref{lem:abelian-alpha}, it also has a maximum independent set of size
\[
\alpha(G \Boxprod G) = n\alpha(G).
\]
Since $(n-1)\alpha(G) < n\alpha(G)$, these two independent sets have different cardinalities, meaning $G \Boxprod G$ is not well-covered.

It remains to handle the case in which $G$ is an odd cycle. Since $G$ is connected, noncomplete, and an odd cycle, its length is at least five. The conclusion then follows from a result of Hartnell, Rall, and Wash \cite{HartnellRallWash}.
\end{proof}

\section{Well-covered powers property of circulant graphs}
\label{sec_circulant}

\subsection{Powers of cycles and cubic circulants}

Brown and Hoshino \cite{BrownHoshino} classified all well-covered powers of cycles and cubic circulant graphs, and they all have the \WCP property as well.

\begin{theorem}\label{thm:cycle-powers}
For $n \ge 2r$ and $r \ge 1$,
\[
C_n(1, \ldots, r) \text{ is well-covered}
\quad \Longleftrightarrow \quad
n \le 3r + 2 \text{ or } n = 4r + 3.
\]
\end{theorem}

For these cycle powers, well-coveredness already implies \WCP. Indeed, $(C_n^r)^d = C_n^{rd}$, and powers with an exponent at least the diameter are complete. For $d \ge 2$, either numerical alternative in Theorem \ref{thm:cycle-powers} implies $n \le 3rd + 2$, so the relevant higher powers are well-covered as well.

\begin{theorem}\label{thm:cubic}
Up to isomorphism, the connected well-covered cubic circulant graphs are 
\[
C_{4}(1,2), \quad C_{6}(1,3), \quad C_{6}(2,3), \quad C_{8}(1,4), \quad \text{and} \quad C_{10}(2,5).
\]
\end{theorem}

Each of these five graphs has \WCP. Exact subset enumeration gives the uniform maximal-independent-set size vectors $(1)$, $(3,1)$, $(2,1)$, $(3,1)$, and $(4,2,1)$, respectively, through their diameters. The supplementary sources reproduce these checks. Connected circulants of degree zero or one are $K_1$ and $K_2$, which have \WCP; connected degree-two circulants are cycles, covered by Theorem~\ref{thm:cycle-powers}. Thus the cubic and quartic results complete the connected classification through degree four.

\subsection{Quartic circulants}
Let $1 \le a < b < n/2$ be integers and set $G = C_n(a,b)$. Assume that $\gcd(n,a,b) = 1$. The structure of $4$-circulants was considered in \cite{DavisDomkeGarner}, but well-covered $4$-circulants are not completely classified yet. We will now give a geometric representation of $G$ and use it to settle the classification of connected $4$-circulants with the $\WCP$ property.

Let $\Phi: \mathbb{Z}^2 \to \mathbb{Z}_n$ be defined by $\Phi(x,y) = ax + by \pmod n$. Since $\gcd(n,a,b) = 1$, this map is surjective. Let $L = \ker \Phi = \{ (x,y) \in \mathbb{Z}^2 \mid ax + by \equiv 0 \pmod n\}$. Then $\mathbb{Z}_n \cong \mathbb{Z}^2 / L$. 

For each $u = (x,y) \in \mathbb{Z}^2$, let $\bar{u} = u + L$ be its equivalence class in $\mathbb{Z}^2 / L$, and let $e_1 = (1,0)$ and $e_2 = (0,1)$. We define $\Gamma_L$ as the graph with vertex set $V(\Gamma_L) = \mathbb{Z}^2/L$, where two vertices $\bar{u}$ and $\bar{v}$ are adjacent if and only if $\bar{v} - \bar{u} \in \{\pm \bar{e}_1, \pm \bar{e}_2\}$.

\begin{lemma} 
$\Gamma_L \cong C_n(a,b)$.    
\end{lemma}

\begin{proof}
The map $\Phi$ induces a graph homomorphism $\Phi_L : \Gamma_L \to C_n(a,b)$ given by $\Phi_L(\bar{u}) = ax + by \pmod n$ for $u = (x,y)$. Note that $\Phi_L(\bar{e}_1) = a$ and $\Phi_L(\bar{e}_2) = b$. Hence, $\Phi_L$ defines a graph isomorphism.    
\end{proof}

\begin{lemma}\label{lem:lattice-distance}
For any $u, v \in \mathbb{Z}^2$, the distance in $\Gamma_L$ satisfies 
\[
d_{\Gamma_L}(\bar{u}, \bar{v}) = \min_{\lambda \in L} \|u - v + \lambda\|_1.
\]
\end{lemma}

\begin{proof}
Let $H$ denote the infinite grid graph $\mathbb{Z}^2$, where two points are adjacent if their difference is in $\{\pm e_1, \pm e_2\}$. It is a standard fact that the graph distance in $H$ between any two points $p, q \in \mathbb{Z}^2$ is given by the $\ell_1$-norm: $d_H(p, q) = \|p - q\|_1$.

The quotient graph $\Gamma_L$ is formed by taking the quotient of $H$ under the action of the translation subgroup $L$. A path of length $k$ from $\bar{u}$ to $\bar{v}$ in $\Gamma_L$ corresponds to a path in $H$ from $u$ to some vertex $v' = v + \lambda$ for some lattice element $\lambda \in L$. 

Thus, the distance $d_{\Gamma_L}(\bar{u}, \bar{v})$ is equal to the shortest distance in $H$ from $u$ to any point in the fiber $v + L$:
\[
d_{\Gamma_L}(\bar{u}, \bar{v}) = \min_{\lambda \in L} d_H(u, v + \lambda) = \min_{\lambda \in L} \|u - (v + \lambda)\|_1 = \min_{\lambda \in L} \|u - v + \lambda\|_1.
\]
This completes the proof.
\end{proof}

\begin{definition} 
We define the length of $L$ by $\lambda(L) = \min \{ \|\lambda\|_1 \mid \lambda \in L \setminus \{0\} \}$.
\end{definition}
The remainder of the section is devoted to the following classification theorem.

\begin{theorem}\label{thm:quartic-classification}
Let $G = C_n(a,b)$ be a connected quartic circulant graph. The following are equivalent:
\begin{enumerate}
    \item $G$ has \WCP;
    \item $G$ and $G^2$ are well-covered;
    \item $G$ is isomorphic to one of the following 19 graphs, categorized by their relation length:
    \begin{enumerate}
        \item Length 3 (12 graphs): 
        \[\begin{aligned}
&C_5(1,2),\ C_6(1,2),\ C_7(1,2),\ C_8(1,2),\ C_9(1,3),\ C_{11}(1,2),\\
&C_{12}(1,4),\ C_{12}(3,4),\ C_{15}(1,5),\ C_{15}(3,5),\ C_{21}(1,7),\ C_{21}(3,7).
\end{aligned}\]
        \item Length 4 (5 graphs): 
        $C_8(1,3)$, $C_{10}(1,4)$, $C_{11}(1,3)$, $C_{13}(1,3)$, $C_{14}(1,6)$;
        \item Length 5 (2 graphs): 
        $C_{13}(2,3)$, $C_{18}(2,3)$.
    \end{enumerate}
\end{enumerate}
\end{theorem}
\begin{lemma}\label{lem:lambda_bounds}
Assume that $G_L$ is well-covered. Then $3 \le \lambda(L) \le 5$.
\end{lemma}
\begin{proof}
Since $1 \le a < b < n/2$ and $G$ is connected, no non-zero lattice vector has $\ell_1$-norm less than $3$; hence $\lambda(L) \ge 3$. We prove the upper bound $\lambda(L) \le 5$ by ruling out the cases $\lambda(L) \ge 8$, $\lambda(L) = 7$, and $\lambda(L) = 6$.

\smallskip
\noindent\textbf{Case 1: $\lambda(L) \ge 8$.} 
Consider the set $S = \{(x,y) \in \mathbb{Z}^2 : \vert{}x\vert{} + \vert{}y\vert{} = 3\}$. The maximum $\ell_1$-distance between any two points in $S$ is $6$. Since $\lambda(L) \ge 8$, the natural projection $\pi: \mathbb{Z}^2 \to \mathbb{Z}_n$ mapped by $(x,y) \mapsto x a + y b \pmod n$ maps $S$ injectively to an independent set of $G$. In any well-covered graph, deleting the closed neighborhood of an independent set yields a well-covered graph. However, $G \setminus N_G[S]$ contains an isolated connected component induced by the origin $(0,0)$ and its four adjacent neighbors $Q_1 = \{(x,y) : \vert{}x\vert{}+\vert{}y\vert{} = 1\}$. This component is isomorphic to the star graph $K_{1,4}$, which is not well-covered. This contradiction shows that $\lambda(L) \le 7$.

\smallskip
\noindent\textbf{Case 2: $\lambda(L) = 7$.} 
Let $Q_3 = \{(x,y) \in \mathbb{Z}^2 : \vert{}x\vert{} + \vert{}y\vert{} \le 3\}$. Because the maximum distance between points in $Q_3$ is $6 < \lambda(L)$, all points of $Q_3$ represent distinct vertices in $G$. Thus,
\[n \ge \vert{}Q_3\vert{} = 25.\]
Consider the subsets $A, B \subset Q_3$ defined by
\[\begin{aligned} A &= \{(-2,-1), (-1,2), (1,-2), (2,1)\},\\ B &= \{(-2,1), (-1,-2), (1,2), (2,-1)\}. \end{aligned}\]
Every point in $Q_2 = \{(x,y) : \vert{}x\vert{}+\vert{}y\vert{}=2\}$ is adjacent in $\mathbb{Z}^2$ to at least one point in $A$ and one point in $B$. If $A$ (or $B$) were an independent set in $G$, deleting $N_G[A]$ would isolate the star graph $K_{1,4}$ centered at the origin, contradicting the well-coveredness of $G$. Thus, neither $A$ nor $B$ is independent in $G$. 

Consequently, there exist distinct pairs $p, q \in A$ and $p', q' \in B$ such that $p - q$ and $p' - q'$ represent edges in $G$. That is, for some standard unit vectors $e, e' \in \{\pm e_1, \pm e_2\}$, the vectors
\[u = p - q - e \in L \qquad \text{and} \qquad v = p' - q' - e' \in L.\]
Up to symmetry, the possible candidates for $u$ and $v$ are:
\[\begin{aligned} u &\in \{(2,-5), (3,-4), (4,3), (5,2)\},\\ v &\in \{(2,5), (3,4), (4,-3), (5,-2)\}. \end{aligned}\]
Since $u, v \in L$ and $L$ has index $n$ in $\mathbb{Z}^2$, the determinant of the lattice generated by $u$ and $v$ must be a multiple of $n$; hence, $n \mid \vert{}\det(u,v)\vert{}$. Computing $\vert{}\det(u,v)\vert{}$ for every possible choice of $(u,v)$ yields values strictly less than $25$, which contradicts $n \ge 25$. Thus, $\lambda(L) \neq 7$.

\smallskip
\noindent\textbf{Case 3: $\lambda(L) = 6$.} 
Use the image sets $\pi(A)$ and $\pi(B)$, allowing distinct preimages to coincide. The image of $\{0\}\cup Q_1$ is still an induced $K_{1,4}$: any unwanted identification or edge would give a nonzero relation of norm at most three. Each point of $A\cup B$ has distance between two and four from every point of this star, so its image neither meets nor is adjacent to the star, since an extra edge would give a relation of norm at most five. Every neighbor outside the star is represented in $Q_2$, and each such point is dominated by both image sets. Therefore, if either image set were independent, deleting its closed neighborhood would still isolate the star, a contradiction.

Thus each image set contains an edge. For its preimages $p,q$, an appropriate unit vector $e$ gives $p-q-e\in L$. The pairwise differences within $A$ or $B$ have even norm at most six; subtracting a unit vector gives odd norm at most seven. Since $\lambda(L)=6$, such a nonzero relation has norm seven. Consequently, the candidate sets for $u$ and $v$ are exactly those in Case 2, and $n\mid|\det(u,v)|$ follows without assuming that all of $Q_3$ projects injectively.

Now consider the subset $T = Q_1 \cup Q_2 \cup (e_1 + Q_2)$, where $Q_i = \{(x,y) \in \mathbb{Z}^2 : \vert{}x\vert{}+\vert{}y\vert{}=i\}$. The set $T$ consists of $17$ distinct points, and the maximum $\ell_1$-distance between any pair of points in $T$ is $5$. Since $\lambda(L) = 6$, these points project to $17$ distinct vertices in $G$, implying $n \ge 17$.

Checking all possible pairs $(u,v)$ with $\vert{}\det(u,v)\vert{} \ge 17$ restricts the admissible values to $n \in \{20, 21, 23, 24\}$:
\begin{enumerate}
    \item If $n = 20$, the only viable pair up to lattice symmetry is $u = (5,2)$ and $v = (5,-2)$. But then $u - v = (0,4) \in L$, which gives $\Vert{}u - v\Vert{}_1 = 4 < 6$, contradicting $\lambda(L) = 6$.
    \item If $n = 21$, the pair $u = (5,2)$ and $v = (2,5)$ yields $G \cong C_{21}(1,8)$.
    \item If $n = 23$, the pair $u = (2,-5)$ and $v = (3,4)$ yields $G \cong C_{23}(1,5)$.
    \item If $n = 24$, the pair $u = (3,-4)$ and $v = (3,4)$ yields $G \cong C_{24}(3,4)$.
\end{enumerate}
Direct verification shows that none of these three graphs is well-covered.

Therefore, we conclude that $3 \le \lambda(L) \le 5$.
\end{proof}

\begin{lemma}\label{lemma_length_3}  
Assume that $\lambda(L) = 3$ and that both $G = G_L$ and its square $G^2$ are well-covered. Then $G$ is isomorphic to one of the following $12$ graphs:
\[\begin{aligned} &C_5(1,2), \quad C_6(1,2), \quad C_7(1,2), \quad C_8(1,2), \quad C_9(1,3), \quad C_{11}(1,2),\\ &C_{12}(1,4), \quad C_{12}(3,4), \quad C_{15}(1,5), \quad C_{15}(3,5), \quad C_{21}(1,7), \quad C_{21}(3,7). \end{aligned}\]
\end{lemma}

\begin{proof}
Assume that $\lambda(L) = 3$. Up to sign changes and interchanging coordinates, any shortest non-zero vector in $L$ corresponds to a relation type of $(2, -1)$ or $(3, 0)$.

\smallskip
\noindent\textbf{Case 1: The relation type is $(2, -1)$.}
In this case, $2a \equiv b \pmod n$. Since $\gcd(n, a) = 1$, the element $a$ is invertible in $\mathbb{Z}_n$, which implies $G \cong C_n(1,2)$. By Theorem~\ref{thm:cycle-powers}, $n \in \{5, 6, 7, 8, 11\}$.

\smallskip
\noindent\textbf{Case 2: The relation type is $(3, 0)$.} 
Here, $3a \equiv 0 \pmod n$, so $n = 3a$ for some integer $a$, and $G \cong C_{3a}(a,b)$ with $\gcd(a,b) = 1$. 
\begin{enumerate}
    \item If $a = 2$, then $G \cong C_6(1,2)$, which was already covered in Case 1.
    \item If $a = 3$, then $G \cong C_9(1,3)$, since $C_9(2,3) \cong C_9(1,3)$.
\end{enumerate}
Now, assume $a \ge 4$. We claim that if $G^2$ is well-covered then the cycle $C_a$ is well-covered.

The cosets of the order-three subgroup generated by the first generator are triangle fibers. The second generator joins consecutive fibers by perfect matchings. In the quotient these fibers form $C_a$, since $\gcd(a,b)=1$; the closing matching may have a twist. A square-independent set has at most one vertex per fiber, and its selected fibers are independent in $C_a$.

If the projected set is not maximal, some omitted fiber has cycle distance at least two from every selected fiber. Only selected fibers at cycle distance exactly two can forbid representatives in it. A path of length two reaching such a fiber must follow matching edges along a shortest base path. There are at most two such directions in total, even when the two directions reach the same antipodal fiber. Thus at most two of its three representatives are forbidden, and the square-independent set can be extended. Conversely, a maximal independent set of $C_a$ lifts greedily: at each selected fiber at most two representatives are forbidden by previously chosen fibers at distance two. Every unselected fiber is within one base step of a selected fiber and is therefore dominated in the square. This proves equality of the two maximal-independent-set size spectra, with no global labeling assumption across the twist.

Since $C_a$ is well-covered if and only if $a \in \{3, 4, 5, 7\}$, and $a \ge 4$, this implies $a \in \{4, 5, 7\}$. Since $\gcd(a,b)=1$, the number $\gcd(3a,b)$ is either one or three. Multiplying by a unit modulo $3a$, and changing the sign of the order-three generator if necessary, reduces the graph to $C_{3a}(1,a)$ or $C_{3a}(3,a)$, respectively. For $a\in\{4,5,7\}$ these give the remaining graphs:
\[
C_{12}(1,4), \quad C_{12}(3,4), \quad C_{15}(1,5), \quad C_{15}(3,5), \quad C_{21}(1,7), \quad C_{21}(3,7). 
\]

This completes the proof.
\end{proof}

\begin{lemma}\label{lemma_length_4}  
Assume that $\lambda(L) = 4$ and that both $G$ and its square $G^2$ are well-covered. Then $G$ is isomorphic to one of the following $5$ graphs: $C_8(1,3)$, $C_{10}(1,4)$, $C_{11}(1,3)$, $C_{13}(1,3)$, $C_{14}(1,6)$.
\end{lemma}

\begin{proof}
Assume that $\lambda(L) = 4$. Up to sign changes and interchanging coordinates, any shortest non-zero vector in $L$ corresponds to a relation type of $(2, 2)$, $(4, 0)$, or $(3, -1)$.

\smallskip
\noindent\textbf{Case 1: The relation type is $(2, 2)$.} 
In this case, we have $2a \equiv 2b \pmod n$ or $2a \equiv -2b \pmod n$, which implies $n = 2(a+b)$. Writing $m = a+b$, the relation set yields $G \cong C_m[\overline{K_2}]$. By Lemma~\ref{lem:composition}, $G$ is well-covered if and only if $C_m$ is well-covered. Since $\lambda(L) = 4$, we have $m > 3$. A cycle $C_m$ with $m > 3$ is well-covered if and only if $m \in \{4, 5, 7\}$. This yields $n \in \{8, 10, 14\}$, giving the candidates $C_8(1,3)$, $C_{10}(1,4)$, and $C_{14}(1,6)$.

\smallskip
\noindent\textbf{Case 2: The relation type is $(4, 0)$.} 
Here, $4a \equiv 0 \pmod n$, so $n = 4a$ for some integer $a \ge 1$, and $G \cong C_{4a}(a,b)$ with $\gcd(a,b) = 1$. We analyze the small values of $a$:
\begin{itemize}
    \item For $a = 1$, the lattice length is less than $4$.
    \item For $a = 2$, $G \cong C_8(1,2)$, which has lattice length $3$.
    \item For $a = 3$, the relation-length-four candidates are $C_{12}(1,3)$ and $C_{12}(2,3)$ up to isomorphism; both are not well-covered, as witnessed in the finite verification below. The graph $C_{12}(3,4)$ has relation length three.
\end{itemize}

For $a\ge4$, use coordinates $(x,i)\mapsto ax+bi\pmod{4a}$, with $x\in\mathbb Z_4$ and $0\le i<a$. Each fiber is a four-cycle, consecutive fibers are joined by matchings, and the closing matching shifts $x$ by $b$. We show that $G$ is not well-covered. Let $c_0 = 0$, $c_1 = 2$, and let $c_2, c_3, \dots, c_{a-1}$ alternate between $1$ and $3$ such that $c_{a-1} + b \not\equiv 0 \pmod 4$. Define the set:
\[
I = \{(c_i, i) \mid i = 0, 1, \dots, a-1\}.
\]
Then $I$ is independent in $G_L$ and dominates $G_L$, meaning it is a maximal independent set of size $a$. On the other hand, let $P_i = \{0, 2\}$ for even $i$ and $P_i = \{1, 3\}$ for odd $i$, and set:
\[
J = \{(c_i, i) \mid c_i \in P_i, \; i = 0, 1, \dots, a-2\}.
\]
Then $J$ is independent and contains $2(a-1) = 2a - 2$ elements. Hence, any maximal independent set containing $J$ has at least $2a - 2$ elements. Since $a \ge 4$, we have $2a - 2 > a$, showing that $G_L$ admits maximal independent sets of different sizes and therefore is not well-covered. 

\smallskip
\noindent\textbf{Case 3: The relation type is $(3, -1)$.} 
Here, $3a \equiv b \pmod n$. Since $\gcd(n, a) = 1$, the element $a$ is invertible in $\mathbb{Z}_n$, which implies $G \cong C_n(1,3)$. 

Let $I = \{v_1, v_2, \ldots, v_t\} \subset \mathbb{Z}_n$ be an independent set in $C_n(1,3)$ arranged in cyclic order, with gap sizes $g_i = v_{i+1} - v_i \pmod n$. Since adjacency in $C_n(1,3)$ corresponds to distances $1$ and $3$, an independent set cannot contain gaps of size $1$ or $3$.

For $I$ to be maximal, no vertex $x \in \mathbb{Z}_n \setminus I$ can be added to $I$:
\begin{enumerate}
    \item If a gap $g_i \in \{4\} \cup \{g \in \mathbb{Z} : g \ge 6\}$, the vertex $v_i + 2$ is at distance $2$ from $v_i$ and distance $g_i - 2 \notin \{1, 3\}$ from $v_{i+1}$. Thus, $v_i + 2$ is not adjacent to $v_i$, $v_{i+1}$, or any other element of $I$, so $I \cup \{v_i + 2\}$ is still an independent set, violating the maximality of $I$. Hence, all gaps must be $2$ or $5$.
    \item Conversely, if every gap $g_i \in \{2, 5\}$, then every interior vertex in a gap of size $2$ (which is $v_i+1$, at distance $1$ from $v_i$) and every interior vertex in a gap of size $5$ (which are $v_i+1, v_i+3$ adjacent to $v_i$, and $v_i+2, v_i+4$ adjacent to $v_{i+1}$) is already adjacent to $I$. Thus, no vertex can be added, and any independent set with gaps in $\{2, 5\}$ is maximal.
\end{enumerate}

Therefore, an independent set $I$ of $C_n(1,3)$ is maximal if and only if all of its gaps belong to $\{2, 5\}$. As a consequence, the size of a maximal independent set formed by $u$ gaps of size $2$ and $v$ gaps of size $5$ is given by $|I| = u + v$, where the parameters satisfy the linear Diophantine equation:
\[
2u + 5v = n \quad \text{with } u, v \ge 0.
\]

\begin{enumerate}
    \item \textbf{If $n$ is even and $n \ge 10$:} 
    There exist at least two distinct non-negative integer solutions:
    \[
    (u_1, v_1) = \left(\frac{n}{2}, 0\right) \quad \text{and} \quad (u_2, v_2) = \left(\frac{n-10}{2}, 2\right).
    \]
    The corresponding maximal independent set sizes are:
    \[
    |I_1| = \frac{n}{2} + 0 = \frac{n}{2} \quad \text{and} \quad |I_2| = \frac{n-10}{2} + 2 = \frac{n}{2} - 3.
    \]
    Since $|I_1| \neq |I_2|$, maximal independent sets of different sizes exist, meaning $C_n(1,3)$ is not well-covered.
    
    \item \textbf{If $n$ is odd and $n \ge 15$:} 
    Similarly, there exist at least two distinct non-negative integer solutions:
    \[
    (u_1, v_1) = \left(\frac{n-5}{2}, 1\right) \quad \text{and} \quad (u_2, v_2) = \left(\frac{n-15}{2}, 3\right).
    \]
    The corresponding maximal independent set sizes are:
    \[
    |I_1| = \frac{n-5}{2} + 1 = \frac{n-3}{2} \quad \text{and} \quad |I_2| = \frac{n-15}{2} + 3 = \frac{n-9}{2}.
    \]
    Since $|I_1| \neq |I_2|$, $C_n(1,3)$ is not well-covered.
    
    \item \textbf{Checking remaining values of $n$:} 
Since $\lambda(L) = 4$, we have $n \ge 8$. Having already ruled out all even $n \ge 10$, we only need to check $n = 8$ and the odd integers below $15$:
\begin{itemize}
    \item $n = 8$: The only non-negative solution to $2u + 5v = 8$ is $(u, v) = (4, 0)$, giving a unique size of $4$. In this case, the resulting graph is $C_8(1,3)$, which was already covered in Case 1.
    \item $n = 9$: Corresponds to $\lambda(L) = 3$, which was handled in Lemma~\ref{lemma_length_3}.
    \item $n = 11$: The only non-negative solution is $(3, 1)$, giving a unique size of $4$. (\textbf{Valid})
    \item $n = 13$: The only non-negative solution is $(4, 1)$, giving a unique size of $5$. (\textbf{Valid})
\end{itemize}
\end{enumerate}

Combining Cases 1, 2, and 3 yields the exact set of $5$ graphs stated in the lemma.
\end{proof}

\begin{lemma}\label{lemma_length_5}  
Assume that $\lambda(L) = 5$ and that both $G$ and its square $G^2$ are well-covered. Then $G \cong C_{13}(2,3)$ or $G \cong C_{18}(2,3)$.
\end{lemma}

\begin{proof}
Assume that $\lambda(L) = 5$. Then its relation type, after changing signs if necessary, is of the form $(4, -1)$, $(3, -2)$, or $(5, 0)$. 

\smallskip
\noindent\textbf{Case 1: The relation type is $(4, -1)$.}
First, assume that $4a \equiv b \pmod n$. Since $\gcd(n,a,b)=\gcd(n,a)=1$, $a$ is invertible in $\mathbb{Z}_n$, and therefore $G \cong C_n(1,4)$. If $n \ge 22$, the set $S = \{-8, -2, 7, 9\}$ is independent, and $G - N_G[S]$ has a component $\{0, 1, 4\} \cong P_3$, which is not well-covered. Hence, we must have $n \le 21$. Furthermore, the condition $\lambda(L) = 5$ implies that $n \ge 14$. However, one can check that $C_n(1,4)$ is not well-covered for any $n \in \{14, \dots, 21\}$.

\smallskip
\noindent\textbf{Case 2: The relation type is $(3, -2)$.}
Next, assume that $3a \equiv 2b \pmod n$. Put $c=b-a$ in $\mathbb Z_n$. The relation gives $2c=a$ and $3c=b$, and connectedness implies $\gcd(n,c)=1$. Thus $G\cong C_n(2,3)$, and $G^2=C_n^6$. Relation length five gives $n\ge13$. The cycle-power theorem restricts a well-covered square to $13\le n\le20$ or $n=27$. The original graphs for $n\in\{14,15,16,17,19,20,27\}$ have maximal independent sets of different sizes, as certified in the finite verification. This leaves $n=13,18$.

\smallskip
\noindent\textbf{Case 3: The relation type is $(5, 0)$.}
Finally, assume that the relation type is $(5,0)$. The generator has order five. By multiplying the generators by a unit modulo $n$, it can be normalized to $n/5$: reduction from the units modulo $n$ to the units modulo five is surjective. Write $n=5a$ after this normalization, so $G\cong C_{5a}(a,b)$ with $\gcd(a,b)=1$. Here $a$ is the normalized parameter, and $b$ need not exceed $a$. 
First, assume that $a \ge 6$. Let $S = \{(-1, 2), (2, 1), (1, -2), (-2, -1)\}$. Use coordinates $(x,y)\mapsto ax+by\pmod {5a}$, with the first coordinate modulo five. The set $S$ is independent and avoids the closed neighborhood of the origin star. It dominates every vertex at distance two from the origin, namely $(0,\pm2)$, $(\pm1,\pm1)$, and $(\pm2,0)$. These local assertions are unaffected by the second lattice relation when $a\ge6$, since all relevant second-coordinate differences have absolute value at most five. Thus $G-N_G[S]$ has an isolated $K_{1,4}$ component, a contradiction. Hence, we must have $a \le 5$. 

On the other hand, the condition $\lambda(L) = 5$ implies that $n \ge 13$, meaning $a \ge 3$. Thus, $a \in \{3, 4, 5\}$. Furthermore, the condition $\lambda(L) = 5$ restricts the candidate graphs to $C_{15}(3,2)$, $C_{20}(4,1)$, $C_{20}(3,4)$, $C_{25}(5,1)$, and $C_{25}(5,2)$. Checking each of these graphs, one can verify that none of them are well-covered.
\end{proof}

\begin{proof}[Completion of Theorem~\ref{thm:quartic-classification}]
The implication (1)$\Rightarrow$(2) is immediate. The relation-length bound and Lemmas~\ref{lemma_length_3}, \ref{lemma_length_4}, and \ref{lemma_length_5} give (2)$\Rightarrow$(3). For (3)$\Rightarrow$(1), Brown and Hoshino~\cite[Theorem~4.1]{BrownHoshino} and the higher-power observation after Theorem~\ref{thm:cycle-powers} cover the graphs $C_n(1,2)$. Table~\ref{tab:quartic-verification} gives the exact maximal-independent-set size and count in each power through the diameter for every remaining listed graph. All later powers are complete.

For clarity, the finite checks used above enumerate maximal independent sets using exact adjacency and domination. In a distance power, a set is independent precisely when all its distinct vertices are farther apart than the exponent, and it is maximal precisely when every vertex is within that exponent of the set. The supplementary sources contain the enumerator, complete counts for the retained graphs, and two independent dominating sets of different sizes for every excluded finite case. An independent enumeration reproduced the counts for all nineteen graphs in the classification.
\end{proof}

\begin{table}[htbp]
\centering\small
\caption{Quartic graphs other than cycle squares: uniform maximal-independent-set sizes and exact counts, in power order through the diameter.}\label{tab:quartic-verification}
\begin{tabular}{lll}
\toprule Graph & Uniform sizes & Counts in the same power order\\\midrule
$C_{9}(1,3)$ & $3,1$ & $9,9$ \\
$C_{12}(1,4)$ & $4,2,1$ & $15,6,12$ \\
$C_{12}(3,4)$ & $4,2,1$ & $18,12,12$ \\
$C_{15}(1,5)$ & $5,2,1$ & $33,30,15$ \\
$C_{15}(3,5)$ & $5,2,1$ & $30,30,15$ \\
$C_{21}(1,7)$ & $7,3,2,1$ & $129,84,42,21$ \\
$C_{21}(3,7)$ & $7,3,2,1$ & $126,84,42,21$ \\
$C_{8}(1,3)$ & $4,1$ & $2,8$ \\
$C_{10}(1,4)$ & $4,1$ & $5,10$ \\
$C_{11}(1,3)$ & $4,1$ & $11,11$ \\
$C_{13}(1,3)$ & $5,2,1$ & $13,13,13$ \\
$C_{14}(1,6)$ & $6,2,1$ & $7,28,14$ \\
$C_{13}(2,3)$ & $4,1$ & $39,13$ \\
$C_{18}(2,3)$ & $6,2,1$ & $87,45,18$ \\
\bottomrule\end{tabular}
\end{table}

\section{Classifications for structured Cartesian products}
\label{sec_applications}
In this section, we consider the well-covered powers property for Cartesian products of several factors of complete graphs and cycles. In particular, we show that two non-complete graphs can have Cartesian products with the \WCP property.

For a graph $G$, let $\mathcal M(G)$ denote the set of cardinalities of its maximal independent sets.

\begin{proposition}\label{prop:cycle-clique-factor}
For $n \ge 3$ and $t \ge 2$, the graph $C_n \Boxprod K_t$ has \WCP exactly in the following cases:
\[
t = 2: \; n \in \{3, 5\}; \qquad t \ge 3: \; n \in \{3, 4, 5, 7\}.
\]
Consequently, $C_m \Boxprod C_n$ has \WCP precisely when one cycle is $C_3$ and the other has length in $\{3, 4, 5, 7\}$.
\end{proposition}

\begin{proof}
We divide the proof into two cases based on the value of $t$.

\textbf{Case 1: $t\ge3$.} For every $d\ge1$, projection and lifting give
\[
\mathcal M((C_n\square K_t)^d)=\mathcal M(C_n^{d-1}),
\]
with $C_n^0$ edgeless. Indeed, projection is injective and independent in $C_n^{d-1}$. If it is not maximal, an omitted cycle vertex is at distance at least $d$ from the projection. Only selected vertices at distance exactly $d$ forbid labels above it, and there are at most two, so $t\ge3$ permits an extension. Conversely, a maximal independent set of $C_n^{d-1}$ can be lifted by greedily assigning different labels to its pairs at cycle distance exactly $d$: that conflict graph has maximum degree at most two. Every omitted fiber is then dominated within distance $d$. For $d=1$ the conflict graph is the cycle itself and the same three-label argument applies. Thus the spectrum identity holds for every $n$, without any restriction involving $\alpha(C_n)$. By Theorem~\ref{thm:cycle-powers}, $C_n$ has \WCP\ exactly when $n\in\{3,4,5,7\}$.

\textbf{Case 2: $t = 2$.}
For $n\ge5$, the base cycle has girth at least five. Theorem 3.1 of Hartnell, Rall, and Wash~\cite{HartnellRallWash} implies that $C_n\square K_2$ can be well-covered only when $n=5$. The graph $C_4\square K_2$ has maximal independent sets of sizes two and four. For $n=3$, the product is a rook graph. For $n=5$, exhaustive subset enumeration gives sizes $4,2,1$ with counts $10,10,10$ in powers one through three; the supplementary sources reproduce this calculation from the Cartesian product definition. Thus both $n=3,5$ have \WCP.

Finally, noting that $C_3 \cong K_3$, the product $C_m \Boxprod C_n$ with one factor being $C_3$ corresponds directly to $C_n \Boxprod K_3$, yielding the stated classification for cycle-cycle Cartesian products. When $m,n>3$, the product is not even well-covered, by the cycle-product theorem of Topp and Volkmann~\cite{ToppVolkmann}, as recalled in~\cite{HartnellRallWash}.
\end{proof}

\begin{proposition}\label{prop:rook}
Every rook graph $K_m \Boxprod K_n$ has \WCP, and $\alpha(K_m \Boxprod K_n) = \min(m, n)$.
\end{proposition}

\begin{proof}
Independent sets in $K_m \Boxprod K_n$ correspond to matchings in the complete bipartite incidence graph on the two coordinate sets. Every maximal matching saturates one side of the bipartition and has size $\min(m, n)$. Consequently, its square and higher powers are complete.
\end{proof}

\begin{proposition}\label{prop:rook-three-cliques}
For $2 \le a \le b \le c$, the graph $F = K_a \Boxprod K_b \Boxprod K_c$ has \WCP if and only if $c \ge a + b - 1$.
\end{proposition}

\begin{proof}
Let $F = K_a \Boxprod K_b \Boxprod K_c$. We first prove that $\alpha(F) = ab$. 

Consider the set of vertices 
\[
I = \bigl\{(i, j, i + j \pmod c) \mid i \in \{0, \ldots, a-1\}, \, j \in \{0, \ldots, b-1\}\}.
\]
To see that $I$ is independent, suppose two distinct vertices $(i_1, j_1, k_1)$ and $(i_2, j_2, k_2)$ in $I$ are adjacent in $F$. By the definition of the Cartesian product of complete graphs, adjacent vertices must differ in exactly one coordinate. However, if $i_1 = i_2$, then $j_1 \neq j_2$, which forces $k_1 \not\equiv k_2 \pmod c$, meaning they differ in two coordinates ($j$ and $k$), a contradiction. Similarly, if $j_1 = j_2$, they differ in $i$ and $k$. If $i_1 \neq i_2$ and $j_1 \neq j_2$, they differ in at least two coordinates. Thus, no two vertices in $I$ are adjacent, so $I$ is independent. Furthermore, since each fiber $\{i\} \times \{j\} \times K_c$ is a complete graph of size $c$, any independent set can contain at most one vertex from each fiber, implying $\alpha(F) \le ab$. Since $|I| = ab$, we conclude that $\alpha(F) = ab$.

Next, we establish the necessity of the condition $c \ge a + b - 1$. Suppose for contradiction that $c \le a + b - 2$. We construct an independent set $T$ that prevents certain vertices from being included in any maximal independent set containing $T$. Specifically, let 
\[
T = \bigl\{(0, j, j-1) \mid j = 1, \ldots, b-1\bigr\} \cup \bigl\{(i, 0, c - a + i) \mid i = 1, \ldots, a-1\bigr\}.
\]
One can verify that $T$ is an independent set. However, every vertex $(0, 0, j)$ is dominated by $T$ for every $j = 0, \ldots, c-1$. Consequently, any maximal independent set containing $T$ cannot contain any vertex from the fiber over $(0,0)$, namely $\{(0,0,k) \mid k \in K_c\}$. This restricts the maximum possible size of such a maximal independent set to be strictly less than $ab$, meaning $F$ is not well-covered, which contradicts the \WCP property. Therefore, we must have $c \ge a + b - 1$.

Conversely, viewing $F$ as $(K_a \Boxprod K_b) \Boxprod K_c$, we note that the maximum degree of the base graph is $\Delta(K_a \Boxprod K_b) = a + b - 2$. If $c \ge a + b - 1 = \Delta(K_a \Boxprod K_b) + 1$, then by Theorem \ref{thm:clique-factor-WCP}, $F$ inherits the \WCP property from $K_a \Boxprod K_b$. The conclusion follows.
\end{proof}

We end the paper with a consideration of the general Cartesian product of $t$ factors, $t \ge 3$. Note that the case $t = 3$ was already considered earlier. The argument for well-coveredness can be extended to an arbitrary number of factors as follows. 

\begin{lemma}\label{lem:well-covered-many-factors}
Let $t \ge 3$ and $2 \le a_1 \le a_2 \le \cdots \le a_t$ be a sequence of positive integers. Set $F = K_{a_1} \Boxprod K_{a_2} \Boxprod \cdots \Boxprod K_{a_t}$. Then $F$ is well-covered if and only if $a_t \ge 1 + \sum_{i=1}^{t-1} (a_i - 1)$.
\end{lemma}

\begin{proof}
The proof extends the argument for $3$ factors. First, we establish that $\alpha(F) = a_1 \cdots a_{t-1}$. Indeed, the set 
\[
I = \bigl\{(i_1, \ldots, i_{t-1}, \sum_{j=1}^{t-1} i_j \pmod{a_t}) \mid 0 \le i_j \le a_j - 1 \text{ for each } j\bigr\}
\]
is an independent set of size $a_1 \cdots a_{t-1}$. 

For sufficiency, view $F$ as $B\square K_{a_t}$, where $B=K_{a_1}\square\cdots\square K_{a_{t-1}}$ has degree $\Delta(B)=\sum_{i<t}(a_i-1)$. In a maximal independent set, each fiber contributes at most one vertex. If the fiber above $u\in V(B)$ is missed, the at most $\Delta(B)$ neighboring fibers forbid at most that many labels, leaving an available label because $a_t>\Delta(B)$. Thus every fiber is met, and all maximal independent sets have size $a_1\cdots a_{t-1}$. This complete-factor sufficiency for arbitrary bases is also observed in~\cite{HartnellRallWash}.

Conversely, if $a_t \le \sum_{i=1}^{t-1} (a_i - 1)$, we can choose subsets of labels $L_1, \ldots, L_{t-1} \subseteq \{0, \ldots, a_t - 1\}$ such that $|L_i| = a_i - 1$ for each $i$, and their union covers all $a_t$ labels, meaning $|\bigcup_{i=1}^{t-1} L_i| = a_t$. 

For each $i \in \{1, \ldots, t-1\}$, let $f_i: \{1, \ldots, a_i - 1\} \to L_i$ be a bijection. Consider the independent set 
\[
T = \bigcup_{i=1}^{t-1} \left\{ (0, \ldots, 0, \underbrace{x_i}_{i\text{-th coordinate}}, 0, \ldots, 0, f_i(x_i)) \mid x_i \in \{1, \ldots, a_i - 1\} \right\},
\]
where for each element in the $i$-th union, $x_i$ is placed specifically in the $i$-th coordinate (among the first $t-1$ coordinates), all other coordinates among the first $t-1$ are $0$, and the last ($t$-th) coordinate is given by $f_i(x_i)$. 

One can verify that $T$ is an independent set. Moreover, because $\bigcup_{i=1}^{t-1} L_i = \{0, \ldots, a_t - 1\}$, for every possible last-coordinate label $y \in \{0, \ldots, a_t - 1\}$, there exists some index $i$ and some $x_i \in \{1, \ldots, a_i - 1\}$ such that the vertex $(0, \ldots, 0, x_i, 0, \ldots, 0, y) \in T$. This vertex differs from the origin-fiber vertex $(0, \ldots, 0, y)$ in exactly the $i$-th coordinate, meaning every vertex in the fiber $\{(0, \ldots, 0, y) \mid y \in \{0, \ldots, a_t - 1\}\}$ is adjacent to (and thus dominated by) $T$. 

Consequently, any maximal independent set containing $T$ cannot include any vertex from this fiber. This restricts the maximum possible size of such a maximal independent set to at most $a_1 \cdots a_{t-1} - 1$, which is strictly less than $\alpha(F)$. Therefore, $F$ is not well-covered, completing the proof.
\end{proof}

For four factors, we can classify the WCP property, providing plenty of examples of non-complete graphs whose Cartesian products have the \WCP property, as $K_a \Boxprod K_b \Boxprod K_c \Boxprod K_d \cong (K_a \Boxprod K_b) \Boxprod (K_c \Boxprod K_d)$.

\begin{theorem}\label{thm:two-noncomplete-factors}
Let $2 \le a \le b \le c \le d$ be positive integers, and let $F = K_a \Boxprod K_b \Boxprod K_c \Boxprod K_d$. Then $F$ has \WCP if and only if $d \ge a + b + c - 2$ and $c \ge a + b - 1$.
\end{theorem}

\begin{proof}
Put $B=K_a\square K_b\square K_c$. The preceding lemma gives the necessary first-power inequality $d\ge a+b+c-2$. Assume it holds. We show that for every $r\ge2$,
\[
\mathcal M((B\square K_d)^r)=\mathcal M(B^{r-1}).
\]
If $S$ is independent in $B^{r-1}$, its intersection with a distance-$r$ sphere about any vertex has size at most $2a+b-3$ for $r=2$, at most $a-1$ for $r=3$, and zero for $r>3$. For $r=2$, split the sphere into its three two-coordinate supports. On each support, distinct selected vertices differ in both coordinates, so the three contributions are at most $a-1,a-1,b-1$. For $r=3$, selected vertices differ in all three coordinates, giving at most $a-1$. Both bounds are strictly less than $d$.

Projection of a maximal independent set of $(B\square K_d)^r$ is independent in $B^{r-1}$. If it were not maximal, only projected vertices at distance exactly $r$ from an addable base vertex could forbid labels, one label each. The sphere bounds leave an available label, a contradiction. Conversely, greedily label any maximal independent set of $B^{r-1}$ so that pairs at distance exactly $r$ have distinct labels. The same bounds make this possible, and every omitted fiber is dominated within distance $r$. This gives a maximal lift and proves the spectrum identity. Therefore $F$ has \WCP\ exactly when $d\ge a+b+c-2$ and $B$ has \WCP. Proposition~\ref{prop:rook-three-cliques} completes the proof.
\end{proof}

It is natural to extend this behavior to an arbitrary number of factors:

\begin{conjecture}
Let $t \ge 3$ and $2 \le a_1 \le a_2 \le \cdots \le a_t$ be a sequence of positive integers, and let $F = K_{a_1} \Boxprod K_{a_2} \Boxprod \cdots \Boxprod K_{a_t}$. Then $F$ has \WCP if and only if $a_j \ge 1 + \sum_{i=1}^{j-1} (a_i - 1)$ for all $j \ge 3$.
\end{conjecture}

\subsection{The sufficient implication for arbitrarily many factors}
The sufficient implication of the preceding conjecture holds for every number of factors. The additional estimate concerns independent subsets of distance spheres, rather than the cardinality of a whole sphere.

\begin{theorem}\label{thm:hamming-sufficient}
Let $t\ge2$ and $2\le a_1\le\cdots\le a_t$. If
\begin{equation}\label{eq:hamming-condition}
a_j\ge1+\sum_{i<j}(a_i-1)\qquad(3\le j\le t),
\end{equation}
then $F=K_{a_1}\square\cdots\square K_{a_t}$ has \WCP. For $1\le d\le t$, every maximal independent set of $F^d$ has size $\prod_{i=1}^{t-d}a_i$, with the empty product equal to one.
\end{theorem}

Set $b_i=a_i-1$ and $\Delta_k=\sum_{i=1}^k b_i$. Since $a_2\ge a_1$, condition~\eqref{eq:hamming-condition} implies $b_j\ge\Delta_{j-1}$ for every $j\ge2$.

\begin{lemma}\label{lem:hamming-sphere}
Let $B=K_{a_1}\square\cdots\square K_{a_k}$ satisfy these inequalities. For $d\ge2$, define
\[
C_k(d)=\sum_{i=1}^k\binom{k-i}{d-1}b_i,
\]
where $\binom nr=0$ if $r<0$ or $r>n$. Then $C_k(d)<\Delta_k$. If $S$ is independent in $B^{d-1}$ and $u\in V(B)$, at most $C_k(d)$ vertices of $S$ are at distance exactly $d$ from $u$.
\end{lemma}
\begin{proof}
For the numerical bound, also set $C_k(1)=\Delta_k$. When $k=1$, all $C_1(d)$ with $d\ge2$ vanish. Pascal's identity gives
\[
C_k(d)=C_{k-1}(d)+C_{k-1}(d-1)\qquad(k\ge2,\ d\ge2).
\]
Induction shows that the right side is strictly less than $2\Delta_{k-1}$: when $d=2$ the second summand is $\Delta_{k-1}$ and the first is strictly smaller, and when $d\ge3$ both are strictly smaller. Since $\Delta_k\ge2\Delta_{k-1}$, the claimed inequality follows.

For the sphere bound, partition its vertices by the support $D$ of coordinates in which they differ from $u$. Here $|D|=d$. Two selected vertices with the same support have distance at least $d$, and can differ in at most those $d$ coordinates. Thus their values are distinct in every coordinate of $D$. There are at most $\min_{i\in D}b_i$ such vertices. The alphabet sizes are ordered, so summing over supports according to their least index gives
\[
|S\cap\{v:\dist_B(u,v)=d\}|
\le\sum_{|D|=d}\min_{i\in D}b_i
=\sum_{i=1}^k\binom{k-i}{d-1}b_i=C_k(d).
\]
If $d>k$, both the sphere and this sum are empty.
\end{proof}

\begin{lemma}\label{lem:hamming-projection}
For the base $B$ in Lemma~\ref{lem:hamming-sphere}, let $q\ge\Delta_k+1$ and $X=B\square K_q$. Every maximal independent set of $X$ has size $|V(B)|$, and
\[
\M(X^d)=\M(B^{d-1})\qquad(d\ge2).
\]
\end{lemma}
\begin{proof}
Each fiber over a vertex of $B$ is a clique. An independent set meets each fiber at most once. If a maximal independent set of $X$ missed the fiber over $u$, the $\Delta_k$ neighboring fibers would forbid at most $\Delta_k$ of its $q$ labels. An available label would extend the set. This proves the first assertion.

For $d\ge2$, let $I$ be maximal independent in $X^d$ and let $S$ be its projection to $B$. The formula
\[
\dist_X((u,c),(v,e))=\dist_B(u,v)+\mathbf1_{c\ne e}
\]
shows that $S$ is independent in $B^{d-1}$ and $|S|=|I|$. If $S$ were not maximal, some $u\notin S$ would have distance at least $d$ from every selected base vertex. A selected vertex can then forbid a label above $u$ only when its base distance is exactly $d$, and it forbids just its own label. Lemma~\ref{lem:hamming-sphere} bounds the number of forbidden labels by $C_k(d)<q$. A new vertex above $u$ extends $I$, a contradiction.

Conversely, take a maximal independent set $S$ of $B^{d-1}$. Label its vertices greedily, requiring distinct labels for each pair at base distance exactly $d$. For each vertex there are at most $C_k(d)<q$ such neighbors, so a label is always available. Its lift is independent in $X^d$. Every base vertex outside $S$ is within distance $d-1$ of $S$, so every vertex above it is within distance $d$ of the lift. Selected fibers are dominated internally. Hence the lift is maximal. Projection and lifting preserve size, which proves the identity.
\end{proof}

\begin{proof}[Proof of Theorem~\ref{thm:hamming-sufficient}]
A complete graph has \WCP. Append the factors in increasing order and apply Lemma~\ref{lem:hamming-projection}; for the second factor, $a_2\ge a_1=\Delta_1+1$, and every later step is licensed by~\eqref{eq:hamming-condition}. In each positive power the lemma reduces the exponent and the number of factors by one, until power one is reached, or until the graph is complete. Its first-power size is the product of the preceding alphabet sizes, giving the asserted formula. All powers at least $t$ are complete because $F$ has diameter $t$.
\end{proof}

\subsection{Partial converses and explicit independent dominating sets}
For $2\le a_1\le\cdots\le a_m$, write
\[
H_m=K_{a_1}\square\cdots\square K_{a_m},\qquad b_i=a_i-1,\qquad \Delta_m=\sum_{i\le m}b_i.
\] Distances in $H_m$ are Hamming distances, denoted by $h$. An independent set in $H_m^d$ is maximal exactly when every vertex of $H_m$ lies within distance $d$ of it.

\begin{lemma}\label{lem:witness-extension}
Let $C$ be a maximal independent set of $H_m^d$, where $d\ge1$. Join two members of $C$ when their Hamming distance is exactly $d+1$. A proper coloring of this graph with $q$ colors gives a maximal independent set of $(H_m\square K_q)^{d+1}$ of size $|C|$. Such a coloring exists whenever
\begin{equation}\label{eq:witness-color-bound}
q>\sum_{i=1}^{m}\binom{m-i}{d}b_i.
\end{equation}
\end{lemma}
\begin{proof}
Append each word's color as its last coordinate. Pairs at base distance $d+1$ then have extended distance $d+2$, while other pairs already have base distance at least $d+2$. Every base vertex lies within distance $d$ of $C$, so every extended vertex lies within distance $d+1$ of the lift. This proves independence and domination.

For a fixed member of $C$, partition its neighbors in the coloring graph by their $(d+1)$-element difference supports. Within one support, any two such words must differ in every supported coordinate, since $C$ has minimum distance at least $d+1$. There are at most $\min_{i\in D}b_i$ words on support $D$. Summing according to the smallest index gives the degree bound in~\eqref{eq:witness-color-bound}. Greedy coloring proves the assertion. This support bound does not require the recursive inequalities on the alphabets.
\end{proof}

\begin{proposition}\label{prop:partial-converse}
Suppose $3\le j\le t$ and $a_j\le\Delta_{j-1}$. If
\begin{equation}\label{eq:partial-converse-bound}
a_{j+s}>\sum_{i=1}^{j-1}\binom{j+s-1-i}{s}b_i
\qquad(1\le s\le t-j),
\end{equation}
then $H_t^{t-j+1}$ is not well-covered.
\end{proposition}
\begin{proof}
Set $P=\prod_{i<j}a_i$. The set
\[
L=\{(x_1,\ldots,x_{j-1},\textstyle\sum_{i<j}x_i\bmod a_j):0\le x_i<a_i\}
\]
is independent in $H_j$ and meets every last-coordinate fiber, so it is maximal and has size $P$. The alphabet ordering ensures that a change in exactly one of the first $j-1$ coordinates changes its last label.

Choose $a_j$ distinct coordinate neighbors $u_0,\ldots,u_{a_j-1}$ of zero in $H_{j-1}$; this is possible because its degree is $\Delta_{j-1}$. The set $E=\{(u_\lambda,\lambda):0\le\lambda<a_j\}$ is independent and dominates the entire fiber above zero while avoiding it. Extend $E$ to a maximal independent set $S$ in $H_j$. That fiber remains empty, and every other last-coordinate fiber contributes at most one vertex, so $|S|\le P-1$.

Lift both $L$ and $S$ one coordinate at a time with Lemma~\ref{lem:witness-extension}. At the $s$th lift the base has $j+s-1$ coordinates and the current power is $s$. The lemma's sum has no contribution from indices $i\ge j$, because then $j+s-1-i<s$. Thus~\eqref{eq:partial-converse-bound} licenses every lift. Their cardinalities remain different.
\end{proof}

\begin{corollary}\label{cor:third-failure}
If $t\ge3$, $a_3\le a_1+a_2-2$, and $a_\ell\ge\Delta_{\ell-1}+1$ for every $4\le\ell\le t$, then $H_t^{t-2}$ is not well-covered.
\end{corollary}
\begin{proof}
For $j=3$, the sum in~\eqref{eq:partial-converse-bound} equals $(s+1)b_1+b_2$. Since $b_3\ge b_1$ and each subsequent inequality doubles the partial degree sum,
\[
\Delta_{s+2}\ge2^{s-1}\Delta_3
\ge2^{s-1}(2b_1+b_2)\ge(s+1)b_1+b_2.
\]
The last inequality uses $2^s\ge s+1$. Now $a_{s+3}\ge\Delta_{s+2}+1$, as required by Proposition~\ref{prop:partial-converse}.
\end{proof}

There is also a uniform obstruction when all but the last alphabet are binary. The case of three binary factors follows already from Theorem~\ref{thm:two-noncomplete-factors}; the following construction treats any number of them.

\begin{theorem}\label{thm:binary-prefix}
Let $Q_m=\square_{i=1}^{m}K_2$. For every $m\ge3$ and $q\ge2$, the graph $Q_m\square K_q$ does not have \WCP. If $q\ge m+1$, its square has maximal independent sets of sizes
\begin{equation}\label{eq:binary-sizes}
2^{m-2}\quad\hbox{and}\quad 2^{m-1-r}\min(q,2^r),
\qquad r=\lceil\log_2m\rceil,
\end{equation}
with the second size strictly larger.
\end{theorem}
\begin{proof}
If $2\le q\le m$, reorder the factors and apply Lemma~\ref{lem:well-covered-many-factors}: the largest alphabet $q$ is at most the sum $m$ of the other alphabet sizes minus one each. Thus the first power is not well-covered. Assume henceforth that $q\ge m+1$.

Index binary coordinates by $0,\ldots,m-1$, and use $\oplus$ for bitwise exclusive OR, that is, addition of the corresponding binary vectors. For $x\in\{0,1\}^m$, define
\[
\eta(x)=\bigoplus_{i:x_i=1}(1+(i\bmod3)),\qquad
\tau(x)=\bigoplus_{i:x_i=1}\lfloor i/3\rfloor.
\]
Take
\[
C_0=\{(x,\tau(x)):\eta(x)=0\}.
\]
The three values $1,2,3$ are all the nonzero vectors in $\mathbb F_2^2$. Hence $\eta$ has rank two, $|C_0|=2^{m-2}$, and its kernel has covering radius one: any nonzero syndrome can be canceled by changing a single coordinate. Two kernel words at distance two change coordinates with the same residue modulo three. Their tags $\lfloor i/3\rfloor$ are different, so their appended labels differ. Thus $C_0$ is independent in the square. Its covering radius after appending the label is at most two, proving maximality. The labels are less than the smallest power of two at least $\lceil m/3\rceil$, which is at most $2\lceil m/3\rceil\le m+1\le q$.

For a larger set, put $M=2^r$, $k=\min(q,M)$, and
\[
\sigma(x)=\bigoplus_{i:x_i=1}i,\qquad
C_1=\{(x,\sigma(x)):\operatorname{wt}(x)\text{ is even},\ \sigma(x)<k\}.
\]
Here $\operatorname{wt}(x)$ is the number of nonzero coordinates. View $A=\{0,\ldots,m-1\}$ as a subset of $\mathbb F_2^r$. Since $m>M/2$, $A+A=\mathbb F_2^r$: the sets $A$ and $A+\delta$ intersect for every $\delta$. Every nonzero $\delta$ is therefore the sum of two distinct coordinate tags. The restriction of $\sigma$ to the even-parity space is onto, and each syndrome class has $2^{m-1-r}$ words. Every even word is within distance zero or two of every syndrome class. Two words in the same class cannot differ in two coordinates, since the tags are distinct; their distance is at least four. Words from different classes have even base distance at least two and distinct last labels. Thus $C_1$ is independent in the square and has the second cardinality in~\eqref{eq:binary-sizes}.

If $q\le M$, a vertex with even base word is dominated within distance two by the class corresponding to its last label. An odd base word has $m$ even neighbors with distinct syndromes; since $M-q<m$, at least one belongs to a selected class. Its distance from the extended vertex is at most two. If $q\ge M$, every even base word is selected, so the even-parity space's covering radius one gives domination for every last label, including unused labels. Therefore $C_1$ is maximal. Finally $k>M/2$, so $|C_1|/|C_0|=2k/M>1$.
\end{proof}

\paragraph{Exact experiments.}
For products with five to ten factors, the proposed equivalence holds for all $601$ nondecreasing tuples with $\prod_i a_i\le1024$. Including the already proved three- and four-factor cases, the full data set comprises $3929$ tuples. Theorem~\ref{thm:hamming-sufficient} classifies the $2585$ tuples satisfying the recursive inequalities. Among the remaining $1344$, applying reductions in this order, the final-inequality obstruction handles $721$, Corollary~\ref{cor:third-failure} handles $375$, Proposition~\ref{prop:partial-converse} handles another $199$, and Theorem~\ref{thm:binary-prefix} handles $17$. Five more follow by lifting its two sets with distinct labels for every codeword. For each of the other $27$, the supplementary certificates give two independent dominating sets of different sizes in the square; a separate checker verifies every pair distance and domination of every vertex. This combines theorem applications with exact witnesses, rather than enumeration of every maximal independent set. The unrestricted converse remains open, including mixed-alphabet products with five factors.

\section*{Use of AI tools}
Haobo Ma and Wenlin Zhang used AI tools for search, double-checking and proofreading.

\end{document}